\documentclass[11pt,letterpaper]{amsart}
\usepackage{amscd}
\usepackage{amssymb}
\usepackage{amsthm}
\usepackage{enumerate}
\usepackage[normalem]{ulem}
\usepackage{mathtools}
\usepackage[usenames,dvipsnames]{xcolor}
\usepackage{stmaryrd}
\usepackage[left=1in,top=1in,right=1in]{geometry}
\usepackage[mathscr]{eucal}
 \usepackage{cite}
\usepackage{upgreek}
\usepackage[bookmarks=false]{hyperref}
\usepackage[T2A]{fontenc}
\usepackage[utf8]{inputenc}
\usepackage{tikz-cd}
\usepackage{mathrsfs}
\usepackage{xcolor}

\newtheorem{thm}{Theorem}

\newtheorem{prob}[thm]{Problem}

\theoremstyle{definition}

\theoremstyle{remark}

\numberwithin{equation}{section}

\begin{document}

\title{50 Open problems: ultraproduct II$_1$ factors}

\author[S. Kunnawalkam Elayavalli]{Srivatsav Kunnawalkam Elayavalli}
\address{\parbox{\linewidth}{Department of Mathematics, University of Maryland, College Park, \\
William E. Kirwan Hall, 4176 Campus Dr, College Park, MD 20742}}
\email{sriva@umd.edu}
\urladdr{https://sites.google.com/view/srivatsavke}

\begin{abstract}
    A collection of 50 open problems around the structure theory of ultraproducts of II$_1$ factors is presented, along with some annotations and references.
\end{abstract}

\maketitle%

Unless otherwise specified, $\mathcal{U}$ will be a non principal ultrafilter on $\mathbb{N}$. For a tracial von Neumann algebra $N$, denote by $N^\mathcal{U}$ the tracial ultrapower.
Denote by $\mathbb{M}^{\mathcal{U}}= \prod_{n\to \mathcal{U}}\mathbb{M}_n(\mathbb{C})$ the tracial matrix ultraproduct with respect to the ultrafilter $\mathcal{U}$, and $\mathbb{M}_\mathcal{U}$ the $C^*$-norm matrix ultraproduct. For a discrete group $G$, denote by $G^{\mathcal{U}}$ the discrete group ultrapower. A pair of II$_1$ factors $N_1$ and $N_2$ are said to be elementarily equivalent if there exists ultrafilters $\mathcal{U}_1$ and $\mathcal{U}_2$ such that $N_1^{\mathcal{U}_1}\cong N_2^{\mathcal{U}_2}$.

\subsection*{Matrix ultraproducts}

\begin{prob}
Classify elementary equivalence classes of $\mathbb{M}^{\mathcal{U}}$ where $\mathcal{U}$ ranges over ultrafilters on $\mathbb{N}$. What about elementary equivalence classes of $\mathbb{M}_{\mathcal{U}}$?
    \end{prob}

The analogue of the above problem for universal sofic groups has been addressed in significant generality in \cite{alekseev2024nonisomorphicuniversalsoficgroups}.

\begin{prob}
    For every Connes embeddable non Gamma II$_1$ factor $N$ does there exist an embedding $\pi: N\hookrightarrow \mathbb{M}^{\mathcal{U}}$ such that $N'\cap \mathbb{M}^{\mathcal{U}}=\mathbb{C}$?\end{prob}
The above ``\emph{ergodic Connes-embedding} problem'' is credited to A. Ioana, and S. Popa. Partial progress includes the case of von Neumann algebras of character rigid property (T) groups with arbitrarily high dimensional irreducible complex unitary representations such as $L(SL_3(\mathbb{Z}))$, and II$_1$ factors $N$ with positive 1-bounded entropy $h(N)>0$, see \cite{jekeltype, Hayes2018}.

\begin{prob}[\textcolor{blue}{Solved by J. Peterson \cite{peterson2026ultraproductapproximationspropertyt}}]
    A non Gamma factor $N$ is said to be pseudocompact if there exists an ultrafilter $\mathcal{U}$ such that $N$ is elementarily equivalent to $\mathbb{M}^{\mathcal{U}}$. Is $L(\mathbb{F}_n)$ pseudocompact for $n\geq 2$? 
\end{prob}

The existence of non psuedocompact non Gamma factors was proved in \cite{exoticCIKE}.

\begin{prob}
    Does there exist a group von Neumann algebra $L(G)$ that is pseudocompact? Does there exist an icc non inner amenable group von Neumann algebra $L(G)$ that is not pseudocompact? (\textcolor{blue}{Second part solved by J. Peterson \cite{peterson2026ultraproductapproximationspropertyt}})
\end{prob}

\begin{prob}
    Let $N$ be a finite index subfactor of the matrix ultraproduct $\mathbb{M}^\mathcal{U}$. Then is  $N$ also isomorphic to a matrix ultraproduct?
\end{prob}

\begin{prob}
Denote by $\mathbb{S}^{\mathcal U}=\prod_{i\to \mathcal U}\mathbb{S}_n$, the ultraproduct of the symmetric groups with the normalized Hamming distance. Clearly, $\mathbb{S}^{\mathcal U}\subset \mathbb{M}^\mathcal{U}$ in a trace preserving way (where the trace on $\mathbb{S}^{\mathcal U}$ is given by one minus the ultralimit of the Hamming distance to the identity elements). Denote by $W^*(\mathbb{S}^{\mathcal U})$ the von Neumann subalgebra of $\mathbb M^\mathcal{U}$ generated by $\mathbb{S}^{\mathcal U}$. Does $W^*(\mathbb{S}^{\mathcal U})$ contain a copy of a matrix ultraproduct? 
\end{prob}

Note that $W^*(\mathbb{S}^{\mathcal U})$ is never isomorphic to a matrix ultraproduct, indeed because the former has a Cartan subalgebra while the latter does not.

\subsection*{Existential embeddings}

Recall that an embedding of $\iota: N\hookrightarrow M$ is existential if there is an embedding of $M\subset N^{\mathcal{U}}$ that restricts to the diagonal embedding under $\iota$. One can define this notion similarly for inclusions of groups.

\begin{prob}
Does there exist a separable $N$ admitting an existential embedding $L(\mathbb{F}_2)\subset N$ and moreover satisfies $N$ is not isomorphic to a non trivial free product? 
\end{prob}

For some mild evidence for the above in the negative, see \cite{jekel2024upgradedfreeindependencephenomena} wherein a general result on free complementation in $L(\mathbb{F}_2)^\mathcal{U}$ is proved, generalizing previous forms of absorption results in the literature, for instance \cite{Popa1983, HoudayerIoana2024, CyrilAOP, FreePinsker}.

\begin{prob}
Does there exist a nonamenable strongly solid $N$ and $M$ containing an existential copy of $N$ such that $M$ admits a Cartan subalgebra? 
\end{prob}

 It is known that if $G$ contains an existential copy of $\mathbb{F}_2$ then it is strongly solid by \cite{ChifanSinclair}, and also satisfies $h(L(G))=\infty$ by \cite{Sridiagonal} (where $h$ denotes Hayes' 1-bounded entropy \cite{Hayes2018}). The problem below is credited to D. Voiculescu for surface groups.
\begin{prob}
    Does there exist a group $G$ that contains an existential inclusion of $\mathbb{F}_2$ such that $L(G)\ncong L(\mathbb{F}_n)$ for some $n\in \mathbb{N}$? Specifically, does there exist a non elementary surface group $G$ such that $L(G)\ncong L(\mathbb{F}_n)$?
\end{prob}

\begin{prob}
    Let $M$ be a separable II$_1$ factor that contains an existential copy of $N$ where $N$ is a nonamenable biexact II$_1$ factor. Then is $M$ prime?
\end{prob}
For the definition of bi-exactness see the work \cite{ding2023biexactvonneumannalgebras}. The above behavior was observed in $L(\mathbb{F}_2)$, see \cite{Sridiagonal}.

\begin{prob}
    Let $N$ be a non amenable II$_1$ factor with Haagerup property, and $N\subset M$ existentially. Then can $M$ have property (T)? 
\end{prob}

The group analogue of the above question has an easy negative answer because of Shalom's theorem. For another interesting conjecture along these lines see J. Peterson's conjecture, below:

\begin{prob}
    Let $N$ be a II$_1$ factor with Haagerup's property, and $B\subset N^\mathcal{U}$, with Property (T). Then are there finite dimensional subalgebras $M_n\subset N$ such that $B\subset \prod_{n\to \mathcal{U}}M_n$?
\end{prob}

Progress and evidence towards the above has been obtained in \cite{peterson2026ultraproductapproximationspropertyt}.

\begin{prob}
    Does there exist an existential inclusion of $L(\mathbb{F}_2)$ into $\mathbb{M}^\mathcal{U}$? Is the random embedding existential, in the sense of Haar-sampling?
\end{prob}

A positive answer to the above question will have applications to rigidity questions concerning free group factors, to be elaborated by us in an upcoming paper. 

Recall that a II$_1$ factor $N$ is said to be existentially closed if whenever it embeds into a II$_1$ factor $M$, it can embed existentially into $M$. Such factors exist by \cite{ECfactors}, also see \cite{IOANA2024110264, GKEPTConj} for concrete constructions.

\begin{prob}
    Does there exist a countable group $G$ such that $L(G)$ is an existentially closed II$_1$ factor? Does there exist an existentially closed II$_1$ factor which is not isomorphic to a group von Neumann algebra? 
\end{prob}

The above question originally appears in \cite{chifan2022tensorproductindecomposabilityresults}, and was kindly communicated to us by I. Chifan. 

\subsection*{McDuff and Property Gamma factors}

\begin{prob}
Does there exist a sequence of separable non Gamma factors $N_i$ such that their ultraproduct $\prod_{i\to \mathcal{U}}N_i$ has Property Gamma?\end{prob}

\begin{prob}[\textcolor{blue}{Solved by J. Peterson \cite{peterson2026ultraproductapproximationspropertyt}}]
Does there exist a pair of Connes embeddable non Gamma factors that are not elementarily equivalent? Does there exist continuum many non Gamma factors that are pairwise not elementarily equivalent?
\end{prob}

There are a pair of non Gamma factors that are not elementarily equivalent \cite{exoticCIKE}. More generally, J. Peterson has found continuum many group factors that are not elementarily equivalent \cite{peterson2026ultraproductapproximationspropertyt}.

\begin{prob}[\textcolor{blue}{Solved by J. Peterson \cite{peterson2026ultraproductapproximationspropertyt}}]
Does there exist an icc non inner amenable group $G$ such that $L(G)$ is not elementarily equivalent to $L(\mathbb{F}_2)$?
\end{prob}

An ideal candidate for the above would be $SL_3(\mathbb{Z})$, and indeed this was shown to be the case by J. Peterson in \cite{peterson2026ultraproductapproximationspropertyt}.

Recall that a II$_1$ factor is super McDuff if $N'\cap N^{\mathcal{U}}$ is a II$_1$ factor.

\begin{prob}[\textcolor{blue}{Solved by J. Peterson \cite{peterson2026ultraproductapproximationspropertyt}}]
Are there two non elementarily equivalent Connes embeddable super McDuff II$_1$ factors $N,M$? Are there continuum many?
\end{prob}

Continuum many such factors were obtained in \cite{peterson2026ultraproductapproximationspropertyt} via a Bass-Serre type rigidity result in ultraproducts. 

\begin{prob}
If $N$ is a super McDuff II$_1$ factor, and $N$ is elementarily equivalent to $M$ then is $M$ super McDuff? In particular, does there exist a super McDuff II$_1$ factor which is not uniformly super McDuff?\end{prob}

See \cite{goldbring2023uniformlysupermcduffii1} for the definition of uniformly super McDuff factors, and a partial result in this direction.

\begin{prob}[\textcolor{blue}{Solved by J. Peterson \cite{peterson2026ultraproductapproximationspropertyt}}]
    Does there exist two factors with property Gamma and are not McDuff, that are not elementarily equivalent? Are there continuum many such examples? 
\end{prob}

The above is complementary to \cite{boutonnetchifanioana}, which established many non elementarily equivalent McDuff factors.

\begin{prob}
Is $L(F)$ super McDuff, where $F$ is the well known Thompson's group? \end{prob}

See \cite{Jolissaint1998} where it was proved that these group von Neumann algebras are McDuff.

\begin{prob}
    If $M\cong N_1\overline{\otimes} N_2$, and $M$ is McDuff, does it imply that either $N_1$ or $N_2$ are McDuff? 
\end{prob}

For a partial result in this direction see \cite{Marrakchi_2018}.

\subsection*{Sequential commutation}

For the next few problems, recall the equivalence relation of sequantial commutation \cite{elayavalli2023sequential} on the haar unitaries in a II$_1$ factor. For $u,v\in \mathcal{H}(N)$ haar unitaries, say that $u\sim_{N^\mathcal{U}} v$ if there exists a finite set $w_1, \cdots, w_n$ of haar unitaries in $N^{\mathcal{U}}$, such that $[u, w_1]=[w_1,w_2]= \cdots = [w_{n-1}, w_n]= [w_n,v]=0$. In the case that there is a unique orbit among haar unitaries, additionally recall the commutation diameter $\mathfrak{O}(M) := \sup_{u,v\in \mathcal{H}(N)} p_N(u,v)$, where $p_N(u,v)$ denotes the length of the shortest path of sequentially commuting haar unitaries in the ultrapower connecting $u$ and $v$.

\begin{prob}
Does there exist a non Gamma factor with unique sequential commutation orbit, and has sequential diameter equal to $n$, $n\in \mathbb{N}$? 
\end{prob}

\begin{prob}
    Does there exist a property (T) II$_1$ factor with unique sequential commutation orbit? 
\end{prob}

\begin{prob}
    Does there exist a II$_1$ factor $N$ with no diffuse relative $(T)$ subalgebras and with unique sequential commutation orbit? Is there such an $N$ with additionally the Haagerup property?
\end{prob}

See the recent paper \cite{gao2026exoticfactorsweaklycoarse}, where this is handled for the case of no $T$ subalgebras (the relative $T$ case is not covered so the problem remains open). 

\begin{prob}
    Does there exist a II$_1$ factor $N$ satisfying $h(N)=0$  but $N$ has continuum many sequential commutation orbits? Secondly does there exist a II$_1$ factor $N$ such that $h(N)=0$, and $N$ is not generated as a von Neumann algebra by a single sequential commutation orbit?
\end{prob}
\begin{prob}
Let $M= N_1*_B N_2$, and $u_i\in \mathcal{H}(N_i)$ such that $E_B(u_i^k)=0$ for all $k\in \mathbb{N}$. Suppose $v_i\in \mathcal{H}(M^\mathcal{U})$ such that $u_i\sim_{M^\mathcal{U}} v_i$ and $E_{B^\mathcal{U}}(v_i^k)=0$ for all $k \in \mathbb{N}$. Then is $v_1$ freely independent from $v_2$?
\end{prob}

See \cite{HoudayerIoana2024} where this is proven for commutants. 

\begin{prob}
    Let $N$ be a non-prime II$_1$ factor. Does it have a unique sequential commutation orbit? 
\end{prob}

\begin{prob}
    Let $N$ be a separable II$_1$ factor. Then do $N$ and $N^t$ have the same number of sequential commutation orbits?
\end{prob}

\subsection*{Lifting problems}

Recall $A,B\subset M$ are $k$-independent if for all alternating centered words $w$ of length $k$ in $A,B$, one has $\tau(w)=0$. Usage of plain independent refers to $2$-independent. Additionally, $u,v\in \mathcal{U}(M)$ are $k$-independent means that the von Neumann algebras they generate are $k$-independent.

\begin{prob}
Suppose $u,v\in N^\mathcal{U}$ independent unitaries, then can they be lifted to sequences of pairwise independent unitaries in $N$?  
\end{prob}

This is done in some particular cases in \cite{exoticCIKE}.

\begin{prob}[\textcolor{blue}{Solved by G. Patchell and A. Shiner, paper in preparation}]
Suppose $u,v\in N^\mathcal{U}$ are freely independent unitaries, then can they be lifted to sequences of pairwise $k$-independent unitaries in $N$ for all $k$?
\end{prob}
It is proved that 2-independent unitaries can lift to independent unitaries in \cite{HoudayerIoana2024}. 

\begin{prob}
Suppose $M$ is a separable subfactor admitting two embeddings $\pi_i:M\to N^\mathcal{U}$ and further suppose $\varphi_n:N\to N$ are a sequence of unital completely positive maps, such that $\Phi= (\varphi_n)_{\mathcal{U}}$ conjugates $\pi_1$ to $\pi_2$. Then does there exist a unitary $u\in N^\mathcal{U}$ that conjugates $\pi_1$ and $\pi_2$?
\end{prob}

\begin{prob}
Is there a non amenable separable II$_1$ factor $N$ such that any pair of embeddings of $N$ into $N^{\mathcal{U}}$ are conjugate by a unitary in $N^{\mathcal{U}}$? \end{prob}

\begin{prob}
If a Connes embeddable separable II$_1$ factor $N$ satisfies that any pair of embeddings of $N$ into $R^{\mathcal{U}}$ are conjugate by an automorphism of $R^{\mathcal{U}}$, then is $N\cong R$? \end{prob}

See \cite{ScottSri2019ultraproduct, AGKE} for a detailed discussion of such Jung-type properties (\cite{JungTubularity}) and what is known and not known.

\begin{prob}
Is the pentagon right angled Artin group II$_1$-factor stable? \end{prob} See Definition 1.3 in \cite{scottpaper} for the definition of stability. The above question is also open in the sense of (flexible) Hilbert-Schmidt stability.

\subsection*{Miscellaneous}

\begin{prob}[II$_1$ factor Tarski's problem]
    Are the factors $L(\mathbb{F}_n)$ elementarily equivalent for different values of $n$?
\end{prob}
The reduced $C^*$-algebra version of this problem has been settled in the negative in the recent work \cite{elayavalli2025negativeresolutioncalgebraictarski}.

The problem below is due to S. Popa:

\begin{prob}
    Let $N$ be a Connes embeddable separable II$_1$ factor. Does it admit an embedding into $R^\mathcal{U}$ with factorial commutant?
\end{prob}

\begin{prob}
    Let $N_0$ be an index $2$ subfactor of an ultrapower $N^\mathcal{U}$. Then is $N_0\cong \prod_{i\to \mathcal{U}} N_i$ where $N_i\subset N$ are index $2$ subfactors?
\end{prob}

\begin{prob}
    Is $N^{\mathcal{U}}$ ever isomorphic to a group von Neumann algebra? 
\end{prob}

It follows from arguments of S. Popa \cite{PopaOrthoPairs} that ultrapowers can never be non trivial tensor products, crossed products, amalgamated free products, or have diffuse regular abelian subalgebras. 

\begin{prob}
Does there exist a II$_1$ factor $N$ such that $h(N^{\mathcal{U}})=\infty$ and $h(N)=0$, where $h$ denotes 1-bounded entropy?\end{prob}

\begin{prob}
If $(N'\cap N^{\mathcal{U}})'\cap N^{\mathcal{U}}=N$, then is $N\cong R$?  
\end{prob}

The above ``\emph{ultrapower-bicommutant} problem'' is credited to S. Popa, see \cite{Popaindep}.

\begin{prob}
    If $M$ is a II$_1$ factor, then does $\mathbb{M}_2(\mathbb{C}) \otimes M$ always embed into $M^\mathcal U$? (Note this is immediate if $M$ is Connes embeddable.)
\end{prob}

\begin{prob}
Does there exist an enforceable II$_1$ factor? 
\end{prob}

Please consult \cite{enforceable} for the definition. The non-existence of such a factor would imply a negative answer to Connes embedding (which has already been solved in the negative \cite{ji2020mipre}). 

Recall that a II$_1$ factor $N$ is said to be solid if any subalgebra $A\subset N$ with diffuse commutant in $N$ is amenable. Recall that a II$_1$ factor $N$ is said to be $\mathcal{U}$-solid if any subalgebra $A\subset N$ with diffuse commutant in $N^\mathcal{U}$ is amenable. 

\begin{prob}
Is solidity equivalent to $\mathcal{U}$-solidity in II$_1$ factors?\end{prob}

See Proposition 7 in \cite{OzawaSolidActa} (see also \cite{HURigid}). The problem is still open.

\begin{prob}
Is there a pair of II$_1$ factors with the same 2-quantifier theory but are not elementarily equivalent? \end{prob}
For the definition of 2-quantifier theory, please consult \cite{quantifiertheory}.

\begin{prob}
Does there exist II$_1$ factors $N,M$ that are elementarily equivalent, but not as as $C^*$-algebras?
\end{prob}

\begin{prob}
Is every seperable II$_1$ factor selfless as a $C^*$-algebra? 
\end{prob}

For a definition of selfless see \cite{robertselfless}. This is known for all ultraproduct II$_1$ factors, and some select other examples. Also recall that it is known that strict comparison holds for all II$_1$ factors viewed as $C^*$-algebras.

\begin{prob}
    Does there exist a II$_1$ factor $N$ such that $N^\mathcal{U}$ has non-full fundamental group? 
\end{prob}

\begin{prob}
    Suppose $N_1, N_2, M_1, M_2$ are separable II$_1$ factors such that $N_1^\mathcal{U}\overline{\otimes} N_2^\mathcal{U}\cong M_1^\mathcal{U}\overline{\otimes} M_2^\mathcal{U}$, then is $N_i$ elementarily equivalent to $M_i$ up to permutation? 
\end{prob}

\begin{prob}
    Suppose $I,J$ are finite sets, and $\{N_i\}_{i\in I}, \{M_j\}_{j\in J}$ are separable II$_1$ factors such that $*_{i\in I}N_i^\mathcal{U}\cong *_{j\in J}M_j^\mathcal{U}$, then is $I\cong J$ and the factors are mutually elementarily equivalent up to permutation?  
\end{prob}

\subsection*{Acknowledgments} We thank I. Chifan, D. Gao, I. Goldbring, A. Ioana, D. Jekel, A. Marrakchi, G. Patchell, J. Peterson, C. Schafhauser for their extremely helpful inputs. We particularly thank D. Gao and G. Patchell for their very detailed feedback, their insights, and also their encouragement.

\bibliographystyle{plain}
\bibliography{comparison}

@misc{ji2020mipre,
      title={MIP*=RE}, 
      author={Zhengfeng Ji and Anand Natarajan and Thomas Vidick and John Wright and Henry Yuen},
      year={2020},
      eprint={2001.04383},
      archivePrefix={arXiv},
      primaryClass={quant-ph}
}

@article{elayavalli2023sequential,
    AUTHOR = {Kunnawalkam Elayavalli, Srivatsav and Patchell, Gregory},
     TITLE = {Sequential commutation in tracial von {N}eumann algebras},
   JOURNAL = {J. Funct. Anal.},
  FJOURNAL = {Journal of Functional Analysis},
    VOLUME = {288},
      YEAR = {2025},
    NUMBER = {4},
     PAGES = {Paper No. 110719, 28},
      ISSN = {0022-1236,1096-0783},
   MRCLASS = {46L10 (46L35)},
  MRNUMBER = {4832103},
       DOI = {10.1016/j.jfa.2024.110719},
       URL = {https://doi.org/10.1016/j.jfa.2024.110719},
}

@article{exoticCIKE,
	author = {Chifan, Ionu{\c t} and Ioana, Adrian and Kunnawalkam Elayavalli, Srivatsav},
	date = {2023/06/17},
	doi = {10.1007/s00039-023-00649-4},
	id = {Chifan2023},
	isbn = {1420-8970},
	journal = {Geometric and Functional Analysis},
	title = {An exotic {II}$_1$ factor without property {G}amma},
	url = {https://doi.org/10.1007/s00039-023-00649-4},
	year = {2023}}

@article {Popaindep,
    AUTHOR = {Popa, Sorin},
     TITLE = {Independence properties in subalgebras of ultraproduct {$\rm
              II_1$} factors},
   JOURNAL = {J. Funct. Anal.},
  FJOURNAL = {Journal of Functional Analysis},
    VOLUME = {266},
      YEAR = {2014},
    NUMBER = {9},
     PAGES = {5818--5846},
      ISSN = {0022-1236,1096-0783},
   MRCLASS = {46L10},
  MRNUMBER = {3182961},
MRREVIEWER = {Paul\ Skoufranis},
       DOI = {10.1016/j.jfa.2014.02.004},
       URL = {https://doi.org/10.1016/j.jfa.2014.02.004},
}

@article {ChifanSinclair,
    AUTHOR = {Chifan, Ionut and Sinclair, Thomas},
     TITLE = {On the structural theory of {${\rm II}_1$} factors of
              negatively curved groups},
   JOURNAL = {Ann. Sci. \'{E}c. Norm. Sup\'{e}r. (4)},
  FJOURNAL = {Annales Scientifiques de l'\'{E}cole Normale Sup\'{e}rieure. Quatri\`eme
              S\'{e}rie},
    VOLUME = {46},
      YEAR = {2013},
    NUMBER = {1},
     PAGES = {1--33 (2013)}
}

@misc{elayavalli2025negativeresolutioncalgebraictarski,
      title={Negative resolution to the {$C^*$}-algebraic Tarski problem}, 
      author={Srivatsav Kunnawalkam Elayavalli and Christopher Schafhauser},
      year={https://arxiv.org/abs/2503.10505},
      eprint={2503.10505},
      archivePrefix={arXiv},
      primaryClass={math.OA},
      url={https://arxiv.org/abs/2503.10505}, 
}

@article {robertselfless,
    AUTHOR = {Robert, Leonel},
     TITLE = {Selfless {${\rm C}^*$}-algebras},
   JOURNAL = {Adv. Math.},
  FJOURNAL = {Advances in Mathematics},
    VOLUME = {478},
      YEAR = {2025},
     PAGES = {Paper No. 110409, 28},
      ISSN = {0001-8708,1090-2082},
   MRCLASS = {46L05 (46L35 46L54)},
  MRNUMBER = {4924062},
       DOI = {10.1016/j.aim.2025.110409},
       URL = {https://doi.org/10.1016/j.aim.2025.110409},
}

@article {Hayes2018,
	author = {Hayes, Ben},
	title = {1-bounded entropy and regularity problems in von {N}eumann algebras},
   	journal = {Int. Math. Res. Not. IMRN},
	fjournal = {International Mathematics Research Notices. IMRN},
	year = {2018},
	number = {1},
	pages = {57--137},
	issn = {1073-7928},
	mrclass = {37A55 (46L40)},
	mrnumber = {3801429},
	doi = {10.1093/imrn/rnw237},
	url = {https://doi.org/10.1093/imrn/rnw237},
}

@article {ScottSri2019ultraproduct,
    AUTHOR = {Atkinson, Scott and Kunnawalkam Elayavalli, Srivatsav},
     TITLE = {On ultraproduct embeddings and amenability for tracial von
              {N}eumann algebras},
   JOURNAL = {Int. Math. Res. Not. IMRN},
  FJOURNAL = {International Mathematics Research Notices. IMRN},
      YEAR = {2021},
    NUMBER = {4},
     PAGES = {2882--2918},
      ISSN = {1073-7928},
   MRCLASS = {46L10 (46L36 46M07)},
  MRNUMBER = {4218341},
       DOI = {10.1093/imrn/rnaa257},
       URL = {https://doi.org/10.1093/imrn/rnaa257},
}

@article {CyrilAOP,
    AUTHOR = {Houdayer, Cyril},
     TITLE = {Gamma stability in free product von {N}eumann algebras},
   JOURNAL = {Comm. Math. Phys.},
  FJOURNAL = {Communications in Mathematical Physics},
    VOLUME = {336},
      YEAR = {2015},
    NUMBER = {2},
     PAGES = {831--851},
      ISSN = {0010-3616},
   MRCLASS = {46L10 (46L09)},
  MRNUMBER = {3322388},
MRREVIEWER = {Ping Zhong},
       DOI = {10.1007/s00220-014-2237-0},
       URL = {https://doi.org/10.1007/s00220-014-2237-0},
}

@article {jekeltype,
    AUTHOR = {Jekel, David},
     TITLE = {Covering entropy for types in tracial {$\rm W^*$}-algebras},
   JOURNAL = {J. Log. Anal.},
  FJOURNAL = {Journal of Logic and Analysis},
    VOLUME = {15},
      YEAR = {2023},
     PAGES = {Paper No. 2, 68},
      ISSN = {1759-9008},
   MRCLASS = {03C66 (46L51 46L54 94A17)},
  MRNUMBER = {4600403},
MRREVIEWER = {Stanis\l aw\ Goldstein},
       DOI = {10.4115/jla.2023.15.2},
       URL = {https://doi.org/10.4115/jla.2023.15.2},
}

@misc{alekseev2024nonisomorphicuniversalsoficgroups,
      title={On non-isomorphic universal sofic groups}, 
      author={Vadim Alekseev and Andreas Thom},
      year={2024},
      eprint={2406.06741},
      archivePrefix={arXiv},
      primaryClass={math.GR},
      url={https://arxiv.org/abs/2406.06741}, 
}

@misc{jekel2024upgradedfreeindependencephenomena,
      title={Upgraded free independence phenomena for random unitaries}, 
      author={David Jekel and Srivatsav Kunnawalkam Elayavalli},
      year={2024, to appear in Trans. Amer. Math. Soc., https://arxiv.org/abs/2404.17114},
      eprint={2404.17114},
      archivePrefix={arXiv},
      primaryClass={math.OA},
      url={https://arxiv.org/abs/2404.17114}, 
}

@article {OzawaSolidActa,
    AUTHOR = {Ozawa, Narutaka},
     TITLE = {Solid von {N}eumann algebras},
   JOURNAL = {Acta Math.},
  FJOURNAL = {Acta Mathematica},
    VOLUME = {192},
      YEAR = {2004},
    NUMBER = {1},
     PAGES = {111--117},
      ISSN = {0001-5962},
   MRCLASS = {46L10},
  MRNUMBER = {2079600},
MRREVIEWER = {Dorin Ervin Dutkay},
       DOI = {10.1007/BF02441087},
       URL = {https://doi.org/10.1007/BF02441087},
}

@article {FreePinsker,
    AUTHOR = {Hayes, Ben and Jekel, David and Nelson, Brent and Sinclair,
              Thomas},
     TITLE = {A random matrix approach to absorption in free products},
   JOURNAL = {Int. Math. Res. Not. IMRN},
  FJOURNAL = {International Mathematics Research Notices. IMRN},
      YEAR = {2021},
    NUMBER = {3},
     PAGES = {1919--1979},
      ISSN = {1073-7928},
   MRCLASS = {46L54 (60B20)},
  MRNUMBER = {4206601},
       DOI = {10.1093/imrn/rnaa191},
       URL = {https://doi.org/10.1093/imrn/rnaa191},
}

@article{AGKE,
title = {Factorial relative commutants and the generalized {J}ung property for {II}$_1$ factors},
journal = {Advances in Mathematics},
volume = {396},
pages = {108107},
year = {2022},
issn = {0001-8708},
doi = {https://doi.org/10.1016/j.aim.2021.108107},
url = {https://www.sciencedirect.com/science/article/pii/S0001870821005466},
author = { Atkinson, Scott and  Goldbring, Isaac and Kunnawalkam Elayavalli, Srivatsav}
}

@article{scottpaper,
 ISSN = {00222518, 19435258},
 URL = {https://www.jstor.org/stable/27114671},
 author = {Scott Atkinson},
 journal = {Indiana University Mathematics Journal},
 number = {3},
 pages = {pp. 1167--1187},
 publisher = {Indiana University Mathematics Department},
 title = {Some Results on Tracial Stability and Graph Products},
 urldate = {2025-11-20},
 volume = {70},
 year = {2021}
}

@article{IOANA2024110264,
title = {Existential closedness and the structure of bimodules of II1 factors},
journal = {Journal of Functional Analysis},
volume = {286},
number = {4},
pages = {110264},
year = {2024},
issn = {0022-1236},
doi = {https://doi.org/10.1016/j.jfa.2023.110264},
url = {https://www.sciencedirect.com/science/article/pii/S0022123623004214},
author = {Adrian Ioana and Hui Tan}
}

@article{Marrakchi_2018, title={Stability of products of equivalence relations}, volume={154}, DOI={10.1112/S0010437X18007388}, number={9}, journal={Compositio Mathematica}, author={Marrakchi, Amine}, year={2018}, pages={2005–2019}}

@article {quantifiertheory,
    AUTHOR = {Goldbring, Isaac and Hart, Bradd},
     TITLE = {Properties expressible in small fragments of the theory of the hyperfinite {II}$_1$ factor},
   JOURNAL = {Confluentes Math.},
  FJOURNAL = {Confluentes Mathematici},
    VOLUME = {12},
      YEAR = {2020},
    NUMBER = {2},
     PAGES = {37--47},
      ISSN = {1793-7434},
   MRCLASS = {03C66 (46L10)},
  MRNUMBER = {4236701},
MRREVIEWER = {Alessandro\ Vignati},
}

@misc{peterson2026ultraproductapproximationspropertyt,
      title={On ultraproduct approximations and property (T) factors}, 
      author={Jesse Peterson},
      year={2026},
      eprint={2605.16669},
      archivePrefix={arXiv},
      primaryClass={math.OA},
      url={https://arxiv.org/abs/2605.16669}, 
}

@misc{chifan2022tensorproductindecomposabilityresults,
      title={Tensor product indecomposability results for existentially closed factors}, 
      author={Ionut Chifan and Daniel Drimbe and Adrian Ioana},
      year={2022},
      eprint={2208.00259},
      archivePrefix={arXiv},
      primaryClass={math.OA},
      url={https://arxiv.org/abs/2208.00259}, 
}

@misc{gao2026exoticfactorsweaklycoarse,
      title={Exotic full factors via weakly coarse bimodules}, 
      author={David Gao and David Jekel and Srivatsav Kunnawalkam Elayavalli and Gregory Patchell},
      year={2026},
      eprint={2602.00930},
      archivePrefix={arXiv},
      primaryClass={math.OA},
      url={https://arxiv.org/abs/2602.00930}, 
}

@article {boutonnetchifanioana,
    AUTHOR = {Boutonnet, R\'emi and Chifan, Ionu\c t{} and Ioana, Adrian},
     TITLE = {I{I{$_1$}} factors with nonisomorphic ultrapowers},
   JOURNAL = {Duke Math. J.},
  FJOURNAL = {Duke Mathematical Journal},
    VOLUME = {166},
      YEAR = {2017},
    NUMBER = {11},
     PAGES = {2023--2051},
      ISSN = {0012-7094,1547-7398},
   MRCLASS = {46L10 (03C20)},
  MRNUMBER = {3694564},
MRREVIEWER = {M.\ J.\ M\polhk aczy\'nski},
       DOI = {10.1215/00127094-0000017X},
       URL = {https://doi.org/10.1215/00127094-0000017X},
}

@article{Jolissaint1998,
author = {Jolissaint, Paul},
journal = {Annales de l'institut Fourier},
language = {eng},
number = {4},
pages = {1093-1106},
publisher = {Association des Annales de l'Institut Fourier},
title = {Central sequences in the factor associated with Thompson’s group $F$},
url = {http://eudml.org/doc/75310},
volume = {48},
year = {1998},
}

@article {JungTubularity,
    AUTHOR = {Jung, Kenley},
     TITLE = {Amenability, tubularity, and embeddings into {$\mathcal{R}^\omega$}},
   JOURNAL = {Math. Ann.},
  FJOURNAL = {Mathematische Annalen},
    VOLUME = {338},
      YEAR = {2007},
    NUMBER = {1},
     PAGES = {241--248},
      ISSN = {0025-5831},
   MRCLASS = {46L10 (46L54)},
  MRNUMBER = {2295511},
       DOI = {10.1007/s00208-006-0074-y},
       URL = {https://doi.org/10.1007/s00208-006-0074-y},
}

@article {PopaOrthoPairs,
    AUTHOR = {Popa, Sorin},
     TITLE = {Orthogonal pairs of {$\ast $}-subalgebras in finite von
              {N}eumann algebras},
   JOURNAL = {J. Operator Theory},
  FJOURNAL = {Journal of Operator Theory},
    VOLUME = {9},
      YEAR = {1983},
    NUMBER = {2},
     PAGES = {253--268},
      ISSN = {0379-4024},
   MRCLASS = {46L10 (16A30 20C07)},
  MRNUMBER = {703810},
MRREVIEWER = {G. A. Elliott},
}

@misc{goldbring2023uniformlysupermcduffii1,
      title={Uniformly Super McDuff II$_1$ Factors}, 
      author={Isaac Goldbring and David Jekel and Srivatsav Kunnawalkam Elayavalli and Jennifer Pi},
      year={2023},
      eprint={2303.02809},
      archivePrefix={arXiv},
      primaryClass={math.OA},
      url={https://arxiv.org/abs/2303.02809}, 
}

@misc{GKEPTConj,
      title={On conjugacy and perturbation of subalgebras}, 
      author={David Gao and Srivatsav Kunnawalkam Elayavalli and Gregory Patchell and Hui Tan},
      year={2024},
      eprint={2403.08072},
      archivePrefix={arXiv},
      primaryClass={math.OA}
}

@article{Popa1983,
    author = "Popa, Sorin",
    title = "Maximal injective subalgebras in factors associated with free groups",
    journal = "Advances in Mathematics",
    volume = "50",
    number = "1",
    pages = "27 -- 48",
    year = "1983",
    issn = "0001-8708",
    doi = "https://doi.org/10.1016/0001-8708(83)90033-6",
    url = "http://www.sciencedirect.com/science/article/pii/0001870883900336"
}

@article {HURigid,
    AUTHOR = {Houdayer, Cyril and Ueda, Yoshimichi},
     TITLE = {Rigidity of free product von {N}eumann algebras},
   JOURNAL = {Compos. Math.},
  FJOURNAL = {Compositio Mathematica},
    VOLUME = {152},
      YEAR = {2016},
    NUMBER = {12},
     PAGES = {2461--2492},
      ISSN = {0010-437X,1570-5846},
   MRCLASS = {46L10 (46L36 46L54)},
  MRNUMBER = {3594283},
MRREVIEWER = {Qihui\ Li},
       DOI = {10.1112/S0010437X16007673},
       URL = {https://doi.org/10.1112/S0010437X16007673},
}

@article {ECfactors,
    AUTHOR = {Farah, I. and Goldbring, I. and Hart, B. and
              Sherman, D.},
     TITLE = {Existentially closed {${\rm II}_1$} factors},
   JOURNAL = {Fund. Math.},
  FJOURNAL = {Fundamenta Mathematicae},
    VOLUME = {233},
      YEAR = {2016},
    NUMBER = {2},
     PAGES = {173--196}
     }

@unpublished{enforceable,
	Author = {Goldbring, Isaac},
	Note = {J. Inst. Math. Jussieu, to appear},
	Title = {Enforceable Operator Algebras},
	Year = {arXiv:1706.09048}}

@article{HoudayerIoana2024, title={Asymptotic freeness in tracial ultraproducts}, volume={12}, DOI={10.1017/fms.2024.93}, journal={Forum of Mathematics, Sigma}, author={Houdayer, Cyril and Ioana, Adrian}, year={2024}, pages={e88}}

@article {Sridiagonal,
    AUTHOR = {Kunnawalkam Elayavalli, Srivatsav},
     TITLE = {Remarks on the diagonal embedding and strong 1-boundedness},
   JOURNAL = {Doc. Math.},
  FJOURNAL = {Documenta Mathematica},
    VOLUME = {28},
      YEAR = {2023},
    NUMBER = {3},
     PAGES = {671--681},
      ISSN = {1431-0635,1431-0643},
   MRCLASS = {46L10 (03C20 20A15)},
  MRNUMBER = {4705597},
       DOI = {10.4171/dm/918},
       URL = {https://doi.org/10.4171/dm/918},
}

@misc{ding2023biexactvonneumannalgebras,
      title={Biexact von {Neumann} algebras}, 
      author={Changying Ding and Jesse Peterson},
      year={2023},
      eprint={2309.10161},
      archivePrefix={arXiv},
      primaryClass={math.OA},
      url={https://arxiv.org/abs/2309.10161}, 
}
\end{document}